\documentclass[journal]{IEEEtran}
\usepackage{cite}

\ifCLASSINFOpdf
\else
\fi
\usepackage{bm}   
\usepackage{amsmath}
\usepackage{amssymb}
\usepackage{indentfirst}  
\usepackage{graphicx}  
\usepackage{subfigure}
\usepackage{multirow} 
\usepackage{booktabs} 
\usepackage{cite}  
\usepackage{threeparttable}   
\usepackage{algorithm}
\usepackage{algpseudocode}
\usepackage{arydshln}
\usepackage{color}
\usepackage{microtype}
\usepackage{cases}
\usepackage{enumerate}

\newtheorem{prop}{Proposition}

\algtext*{EndWhile}
\algtext*{EndIf}
\usepackage{multirow}

\begin{document}
%
\title{{\color{black}Resilient Disaster Recovery Logistics of Distribution
Systems: Co-Optimize Service Restoration with Repair Crew and Mobile Power Source Dispatch}}
%
%
%

\author{Shunbo~Lei,~\IEEEmembership{Member,~IEEE,}
        Chen~Chen,~\IEEEmembership{Member,~IEEE,}
        Yupeng~Li,~\IEEEmembership{Member,~IEEE,}
        and~Yunhe~Hou,~\IEEEmembership{Senior~Member,~IEEE}
\thanks{This work was supported in part by the National Natural Science
Foundation of China under Grant 51677160, and in part by the Research Grant
Council of Hong Kong SAR through the Theme-based
Research Scheme under Project No. T23-701/14-N.
The work of C. Chen was supported by the U.S. Department of Energy (DOE)'s
Office of Electricity Delivery and Energy Reliability.}
\thanks{S. Lei and Y. Hou are with the Department of Electrical and Electronic
Engineering, The University of Hong Kong, Hong Kong; Y. Hou is also with The
University of Hong Kong Shenzhen Institute of Research and Innovation,
Shenzhen 518057, China (email: leishunbo@eee.hku.hk, yhhou@eee.hku.hk).}
\thanks{C. Chen is with the Energy Systems Division, Argonne National
Laboratory, Argonne, IL 60439 USA (email: morningchen@anl.gov).}
\thanks{Y. Li is with the Department of Computer Science, The University of Hong
Kong, Hong Kong (email: ypli@connect.hku.hk).}}

\maketitle
\begin{abstract}
Repair crews(RCs) and mobile power sources(MPSs)
are critical resources for distribution system (DS) outage management
after a natural disaster. However, their logistics is not~well investigated.
{\color{black} We propose a resilient scheme 
for disaster recovery logistics to co-optimize DS restoration with
dispatch of RCs~and MPSs. A novel co-optimization
model is formulated to route~RCs and MPSs in the
transportation network, schedule them in the DS, and reconfigure the DS for
microgrid formation coordinately, etc.}
The model incorporates different timescales of DS~restoration
and RC/MPS dispatch, the coupling of transportation
and power networks, etc. To ensure radiality of
{\color{black} the DS with variable~physical structure and MPS allocation,}
we also model~topology~constraints based on the concept of spanning
forest. The model is convexified equivalently and
linearized into a mixed-integer linear programming. {\color{black}To reduce 
its computation time,} preprocessing methods are proposed
to pre-assign a minimal set of repair tasks to~depots~and reduce
the number of candidate nodes for~MPS~connection.~Resilient
recovery strategies thus are generated to enhance service restoration,
especially by dynamic formation of microgrids that are powered
by MPSs and topologized by repair actions of RCs and network
reconfiguration of the DS. Case studies demonstrate the proposed
methodology.
\end{abstract}


%
\IEEEpeerreviewmaketitle

\vspace{-6.4pt}
\section{Introduction}
%
%
%
%
\IEEEPARstart{R}{ecent}
years have witnessed more frequent natural~disasters causing severe power outages,
which result in~tremendous economic loss, etc.
The urgency to enhance power~grid resilience is highlighted.
Specifically, efficient outage management is one of the critical
requirements for resilient power~grids.

\begin{figure}[t!]
  \centering
  \includegraphics[width=2.4in]{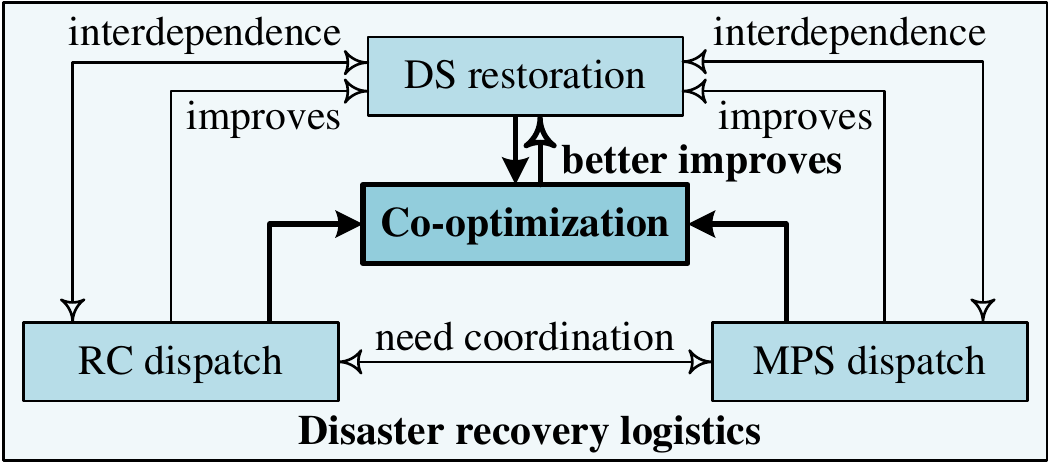}
  \vspace{-2.5mm}
  \caption{Relationships among DS restoration, RC dispatch and MPS dispatch.}
  \vspace{-4pt}
  \label{modules}
\end{figure}

Repair crews (RCs) are crucial response resources for~power grid outage management
against natural disasters.
It is desired that they repair damaged components in an optimal order.~Mobile power sources,
including truck-mounted mobile emergency generators (MEGs) \cite{ZhouBin2017},
mobile energy storage systems~(MESSs) \cite{EPRI2009}, are also important flexibility~resources for grid~restoration.
They can supply critical loads that lose access to the main grid power.
As for power grids, i.e., distribution systems~(DSs) in this work, their restoration involves many~decisions~and~different
strategies, e.g., forming~microgrids~by network reconfiguration \cite{Chen2016}\hspace{0.025cm}\cite{Ding2017}.
Integrating these elements (see Fig. \ref{modules}), DS outage management becomes a
\emph{disaster recovery logistics} problem
to route RCs~and MPSs in the transportation
network \cite{Toth2002},~schedule them in the DS, {\color{black} and reconfigure the DS, etc., for electric service restoration.}
Involving different resources~and infrastructures, this problem is currently not well investigated.

As depicted in Fig. \ref{modules}, the 3 sub-problems are interdependent.
{\color{black}Coordination is required to better improve service restoration.}
References \cite{Perrier2013}\hspace{0.025cm}\cite{Perrier2013_2} 
review the dispatch and planning of RCs~in DSs.
Their interdependence is often simplified or ignored,~e.g., in \cite{Zografos1998}.
Only several publications, e.g., \cite{Ari2017}, co-optimize RC~dispatch with DS restoration.
As for MPSs, in \cite{Abdeltawab2017}, \cite{ChenYY2016} and \cite{GaoHaixiang2017},
MESSs are dispatched to reduce DS operation cost, improve reliability,
and enhance preparedness against natural disasters, respectively.
Review on MEGs can be found in \cite{Lei2016}.
In general, MPS dispatch is seldomly co-optimized with DS restoration, and never in a dynamic manner.
Besides, no publication has considered both RCs and MPSs for DS service restoration.

This paper proposes a co-optimization method for DS~disaster recovery logistics.
RC dispatch and MPS dispatch are jointly coordinated with DS restoration to {\color{black}better}
enhance grid resilience. A non-convex mixed-integer non-linear
programming~(MINLP) model is formulated to co-optimize the routing and scheduling of RCs and MPSs,
dynamic microgrid formation of the DS,~etc.
Issues such as different timescales of RC/MPS dispatch
and DS restoration, the coupling of transportation and power networks, MESSs' state of charge (SoC)
variations over time, and radiality constraints for {\color{black} a DS with variable physical structure and MPS allocation,}
are resolved.
The model is equivalently convexified as a mixed-integer
second-order cone programming (MISCOP) and further linearized to be a mixed-integer
linear programming (MILP). Preprocessing methods, e.g., pre-assigning minimum repair tasks, are also proposed to reduce its computation
time.

In the following,
Section~\ref{RCformulation},~\ref{MPSformulation} and \ref{DSformulation}
formulate RC~dispatch, MPS dispatch and DS restoration, respectively. 
Section \ref{CoOptModel} builds the co-optimization model.
Section \ref{Algorithms}, \ref{cases}~and~\ref{conclusion}
present~the solution method,
case studies and conclusion, respectively.

\section{Routing and Scheduling of RCs}\label{RCformulation}

Dispatch of RCs involves two interdependent sub-tasks, i.e., routing and scheduling.
Routing is to select a route for each RC to travel among depots and damaged components.
Scheduling is to set a timetable for RCs' traveling and repairing behaviors.

Let $\boldsymbol{V}\triangleq\boldsymbol{V_{1}}\cup\boldsymbol{V_{2}}$
be the set of damaged components ($\boldsymbol{V_{1}}$) 
and depots ($\boldsymbol{V_{2}}$).
Let $\boldsymbol{A}\triangleq\{(m,n),\forall m,n\in \boldsymbol{V}\}$
be the set of edges for all pairs of vertices.
Then, the routing of RCs is to find their optimal paths in graph $\boldsymbol{G}\triangleq(\boldsymbol{V},\boldsymbol{A})$.
As the routing problem can be seen as a generalization of the traveling salesman problem (TSP) \cite{Dantzig1959},
its modeling is normally based on the TSP formulation, e.g., in \cite{Ari2017}.
Here we propose a model much simpler and more appropriate for the studied problem.

The TSP formulation is based on edge-wise routing variables $a_{mn}^{k}$ (1 if RC $k$ travels edge $(m,n)$, 0 otherwise).
In this work, we use vertex-wise variables $a_{m,t}^{k}$ instead (1 if RC $k$ is 
at vertex $m$ at time $t$, 0 otherwise).
Note that with subscript $t$, $a_{m,t}^{k}$ are actually scheduling variables.
We will show that the routing problem can be incorporated using $a_{m,t}^{k}$ only.

In scheduling RCs, each RC $k$ can only be visiting at most one vertex at each time $t$:
\begin{equation}\label{AtMostOneRepair}
  \begin{aligned}
    \sum_{m\in\boldsymbol{V}}a_{m,t}^{k}\le1,\forall k,\forall t.
  \end{aligned}
\end{equation}
For conciseness of the objective function to be modeled in Section \ref{CoOptModel},
auxiliary variables $\beta_{t}^{k}$ are introduced:
\begin{equation}\label{AuxRC}
  \begin{aligned}
    \beta_{t}^{k}=1-\sum_{m\in\boldsymbol{V}} a_{m,t}^{k},\forall k, \forall t.
  \end{aligned}
\end{equation}
If RC $k$ is traveling on the transportation network at time $t$, $\beta_{t}^{k}=1$; 
if it is visiting one of the vertices, $\beta_{t}^{k}=0$. 

The routing of RCs imposes travel time constraints on RC scheduling.
For example, if RC $k$ is at vertex $m$ at time $t=1$ (i.e., $a_{m,1}^{k}=1$), and it takes 2 time periods to travel from vertex $m$ to $n$,
then RC $k$ is possible at vertex $n$ only after $t=3$ (i.e., $a_{n,2}^{k}=a_{n,3}^{k}=0$). Such constraints are modeled as:
\begin{equation}\label{TravelTime}
  \begin{aligned}
    a_{n,t+\tau}^{k}+a_{m,t}^{k}\le1,\forall k, \forall m\neq n, \forall \tau\le tr_{mn}, \forall t\le T-\tau.
  \end{aligned}
\end{equation}
where $tr_{mn}$ is the travel time between vertices $m$ and $n$; and $T$ is the number of time periods.
Note that \eqref{TravelTime} actually includes \eqref{AtMostOneRepair} implicitly with $\tau=0$.
To reduce the number of constraints, \eqref{TravelTime} can also be equivalently transformed into:
\begin{equation}\label{TravelTimeTrans}
  \begin{aligned}
    \sum_{{\color{black}\tau=t+1}}^{\min (t+tr_{mn},T)}a_{n,\tau}^{k}\le(1-a_{m,t}^{k})\cdot\min (tr_{mn},T-t),
    \\\forall k, \forall m\neq n, \forall t.
  \end{aligned}
\end{equation}

Actually, \eqref{TravelTime} or \eqref{TravelTimeTrans} is all we need to incorporate the routing sub-problem in the co-optimization.
Other routing constraints, e.g., path-flow balance, will be satisfied implicitly.
The validity of \eqref{TravelTime} and \eqref{TravelTimeTrans} is supported by the following proposition:
\begin{prop}\label{RouteModel}
  For any RC scheduling plan satisfying \eqref{TravelTime} or \eqref{TravelTimeTrans},
  there exists at least one corresponding feasible path in $\boldsymbol{G}$.
\end{prop}

The proof of Proposition \ref{RouteModel} is straightforward.
The key is that we can use $a_{m,t}^{k}$ to retrieve the path of RC $k$ simply by:
\begin{equation}\label{RecoverPath}
  \begin{aligned}
    \begin{bmatrix}
      a_{1,1}^{k} & \cdots & a_{V,1}^{k} \\
      \vdots & \ddots & \vdots \\
      a_{1,T}^{k} & \cdots & a_{V,T}^{k}
    \end{bmatrix}
    \cdot
    \begin{bmatrix}
      1 \\
      \vdots \\
      V
    \end{bmatrix}
  \end{aligned}
\end{equation}
where $V$ is the cardinality of $\boldsymbol{V}$.

{\color{black} The advantages of our proposed model for the routing sub-problem
are mainly twofold: 1) The issues of transportation-DS networks coupling
and their different timescales are resolved in a simpler manner compared with the common TSP formulation.
2) It does not enforce RCs to visit all vertices of damaged components ($\boldsymbol{V_{1}}$).
Thus, in the co-optimization, we can repair a minimal set of damaged components to restore all loads.}

Next, we formulate the repair plan, which depends on the routing and scheduling of RCs.
Let $z_{m,t}^{k}$ be 1 if damaged component $m$ is repaired by RC $k$ at time $t$, 0 otherwise.
Let $rt_{m}^{k}$ be the required time periods for RC $k$ to repair damaged component $m$.
Then, we have:
\begin{equation}\label{Repair1}
  \begin{aligned}
    z_{m,t}^{k}\le\frac{{\color{black} \sum_{\tau=1}^{t}a_{m,\tau}^{k}}}{rt_{m}^{k}},\forall k, \forall m\in\boldsymbol{V_1}, \forall t.
  \end{aligned}
\end{equation}
\begin{equation}\label{Repair1+}
  \begin{aligned}
    z_{m,t}^{k}\le z_{m,t+1}^{k},\forall k, \forall m\in\boldsymbol{V_1}, \forall t\le T-1.
  \end{aligned}
\end{equation}
For example, if it takes $rt_{m}^{k}=2$ time periods of RC $k$ to repair damaged component $m$,
and RC $k$ repairs it at time $t=1\sim2$ (i.e., $[a_{m,1}^{k},a_{m,2}^{k}]=[1,1]$),
then according to \eqref{Repair1} we have $z_{m,1}^{k}\le\frac{1}{2}$ and $z_{m,2}^{k}\le\frac{1+1}{2}$. 
As $z_{m,t}^{k}$ is binary, we further have $[z_{m,1}^{k},z_{m,2}^{k},\cdots,z_{m,T}^{k}]=[0,1,\cdots,1]$.
Note that decisions such as $[a_{m,1}^{k},a_{m,2}^{k},a_{m,3}^{k}]=[1,0,1]$ (idling the RC at $t=2$ during the repair)
and $[a_{m,1}^{k},a_{m,2}^{k},a_{m,3}^{k}]=[1,1,1]$ (idling the RC at $t=3$ after the repair)
are also feasible in \eqref{Repair1}.
However, such decisions are not optimal. They will be eliminated when solving the co-optimization problem.

The following constraints state that each damaged component is repaired only once by one of the RCs:
\begin{equation}\label{Repair2}
  \begin{aligned}
    \sum_{k}z_{m,t}^{k}\le1,\forall m\in\boldsymbol{V_1},\forall t.
  \end{aligned}
\end{equation}
Actually, \eqref{Repair2} is dispensable. Without it, the co-optimization
also seeks an optimal solution without repeated repairs.
Nevertheless, with \eqref{Repair2}, the linear programming relaxation
of our MILP co-optimization model is tightened. It helps to reduce the solution time
using MILP solvers that are generally based on the branch-and-bound framework \cite{Miller2003}.

The following constraints enforce that RC $k$'s resource capacity suffices the total resources required by its repair tasks:
\begin{equation}\label{Repair3}
  \begin{aligned}
    \sum_{m\in\boldsymbol{V_1}}{\color{black}z_{m,T}^{k}}\cdot rs_{m}\le RS^{k},\forall k.
  \end{aligned}
\end{equation}
where $rs_{m}$ is the number of resources required to repair damaged component $m$;
and $RS^{k}$ is RC $k$'s resource capacity.

\section{Routing and Scheduling of MPSs}\label{MPSformulation}

Dispatch of MPSs also involves two interdependent sub-tasks, i.e., routing and scheduling.
Routing is to select a route for each MPS to travel among candidate nodes for MPS connection in the DS. 
Scheduling is to manage MPSs' traveling behaviors and power outputs (or inputs for MESSs when charging) over the considered time window.

Let $\boldsymbol{N}$ be the set of DS nodes;
and $\boldsymbol{N'}\subset\boldsymbol{N}$ be the~set~of~candidate nodes for MPS connection.
Let $\boldsymbol{M_{1}}$ be the set of MEGs;
$\boldsymbol{M_{2}}$ be the set of MESSs;
and $\boldsymbol{M} \triangleq \{{\boldsymbol{M_{1}} \cup \boldsymbol{M_{2}}}\}$
be the set~of MPSs.
Let $vol_{i}$ be the allowed number of MPSs connected to node $i$.
Similar to the RC dispatch model, let $\alpha_{i,t}^{s}$ be 1 if MPS $s$ is connected to node $i$ at time $t$, 0 otherwise.
Then, we have:
\begin{equation}\label{MPS1}
  \begin{aligned}
    \sum_{i\in\boldsymbol{N'}}\alpha_{i,t}^{s}\le1,\forall s,\forall t.
  \end{aligned}
\end{equation}
\begin{equation}\label{MPS2}
  \begin{aligned}
    \sum_{\tau=t}^{\min (t+tr_{ij},T)}\alpha_{j,\tau}^{s}\le(1-\alpha_{i,t}^{s})\cdot\min (tr_{ij},T-t),
    \\\forall s, \forall i\in\boldsymbol{N'}, \forall j\in\boldsymbol{N'}\setminus \{i\}, \forall t.
  \end{aligned}
\end{equation}
\begin{equation}\label{MPS2_3}
  \begin{aligned}
    \sum_{s\in\boldsymbol{M}}\alpha_{i,t}^{s}\le vol_{i},\forall i \in \boldsymbol{N'},\forall t.
  \end{aligned}
\end{equation}
Above, \eqref{MPS1} and \eqref{MPS2} are essentially the same as \eqref{AtMostOneRepair} and \eqref{TravelTimeTrans}, respectively;
and \eqref{MPS2_3} limits the number of MPSs connected to each candidate node by its capacity.
The routing and traveling behaviors of MPSs are sufficiently formulated by \eqref{MPS1}-\eqref{MPS2_3}.
Again, for conciseness of the objective function to be modeled in Section \ref{CoOptModel},
auxiliary variables $\beta_{t}^{s}$ are introduced:
\begin{equation}\label{MPS3}
  \begin{aligned}
    \beta_{t}^{s}=1-\sum_{i\in\boldsymbol{N'}}\alpha_{i,t}^{s},\forall s, \forall t.
  \end{aligned}
\end{equation}
That is, $\beta_{t}^{s}=1$ if MPS $s$ is traveling on the transportation network at time $t$;
$\beta_{t}^{s}=0$ if it is connected to the DS.

Next, we formulate the power dispatch of MPSs.
First, for MEGs, let $gp_{t}^{s}$ and $gq_{t}^{s}$ be the real and reactive power output
of MEG $s$ at time $t$, respectively; let $\overline{gp}^{s}$ and $\overline{gq}^{s}$
be the maximum real and reactive power output of MEG $s$, respectively.
Then, the following constraints are enforced:
\begin{equation}\label{MEG1}
  \begin{aligned}
    0 \le gp_{t}^{s} \le \sum_{i \in \boldsymbol{N'}} \alpha_{i,t}^{s} \cdot \overline{gp}^{s},
    \forall s \in \boldsymbol{M_{1}}, \forall t.
  \end{aligned}
\end{equation}
\begin{equation}\label{MEG2}
  \begin{aligned}
    0 \le gq_{t}^{s} \le \sum_{i \in \boldsymbol{N'}} \alpha_{i,t}^{s} \cdot \overline{gq}^{s},
    \forall s \in \boldsymbol{M_{1}}, \forall t.
  \end{aligned}
\end{equation}
With \eqref{MEG1} and \eqref{MEG2}, MEGs' power outputs are restricted by~their capacities,
and enforced to be zero if a MEG is not connected to the DS.
Fuel limits of MEGs are not included, as they can be refueled by tanker trucks
in case of long-term blackouts~\cite{Iwai2009}.

Second, for MESSs,
let $c_{t}^{s}$ and $d_{t}^{s}$ be 1 if MESS $s$ is charging and discharging at time $t$, respectively, 0 otherwise;
let $cp_{t}^{s}$ and $dp_{t}^{s}$ be the charging and discharging power of MESS $s$ at time $t$, respectively;
let $\overline{cp}^{s}$ and $\overline{dp}^{s}$ be the maximum charging and discharging power of MESS $s$, respectively;
let $soc_{t}^{s}$ be the SoC of MESS $s$ at time $t$;
let $\eta_{s}^{c}$ and $\eta_{s}^{d}$ be the charging and discharging efficiency of MESS $s$, respectively;
let $\underline{soc}^{s}$ and $\overline{soc}^{s}$ be the minimum and maximum SoC of MESS $s$, respectively.
Then, the following constraints are enforced:
\begin{equation}\label{MESS1}
  \begin{aligned}
    c_{t}^{s} + d_{t}^{s} \le \sum_{i \in \boldsymbol{N'}} \alpha_{i,t}^{s},
    \forall s \in \boldsymbol{M_{2}}, \forall t.
  \end{aligned}
\end{equation}
\begin{equation}\label{MESS2}
  \begin{aligned}
    0 \le cp_{t}^{s} \le c_{t}^{s} \cdot \overline{cp}^{s},
    \forall s \in \boldsymbol{M_{2}}, \forall t.
  \end{aligned}
\end{equation}
\begin{equation}\label{MESS3}
  \begin{aligned}
    0 \le dp_{t}^{s} \le d_{t}^{s} \cdot \overline{dp}^{s},
    \forall s \in \boldsymbol{M_{2}}, \forall t.
  \end{aligned}
\end{equation}
\begin{equation}\label{MESS3_4}
  \begin{aligned}
    0 \le gq_{t}^{s} \le (c_{t}^{s} + d_{t}^{s}) \cdot \overline{gq}^{s},
    \forall s \in \boldsymbol{M_{2}}, \forall t.
  \end{aligned}
\end{equation}
\begin{equation}\label{MESS4}
  \begin{aligned}
    soc_{t+1}^{s} = soc_{t}^{s} +(cp_{t}^{s} \cdot \eta_{s}^{c} - \frac{dp_{t}^{s}}{\eta_{s}^{d}}) \cdot \Delta t,
    \\\forall s \in \boldsymbol{M_{2}}, \forall t \le T-1.
  \end{aligned}
\end{equation}
\begin{equation}\label{MESS5}
  \begin{aligned}
    \underline{soc}^{s} \le soc_{t}^{s} \le \overline{soc}^{s},
    \forall s \in \boldsymbol{M_{2}}, \forall t.
  \end{aligned}
\end{equation}
where $\Delta t$ is the duration of one time period.
Specifically, \eqref{MESS1} ensures that in each time period, charging and discharging are mutually exclusive
states of a MESS, and it can neither charge nor discharge if it is not connected to the DS;
\eqref{MESS2} and \eqref{MESS3} specify MESSs' charging and dicharging power limits, respectively
(if a MESS is not in the charging/discharging state, its charging/discharging power is limited to be zero);
\eqref{MESS3_4} is similar to \eqref{MEG2};
\eqref{MESS4} expresses MESSs' SoC variations over time;
and \eqref{MESS5} imposes SoC ranges for MESSs.

\section{Dynamic Network Reconfiguration and \\ Power Dispatch of the DS}\label{DSformulation}

\begin{figure}[t!]
  \centering
  \includegraphics[width=3.1in]{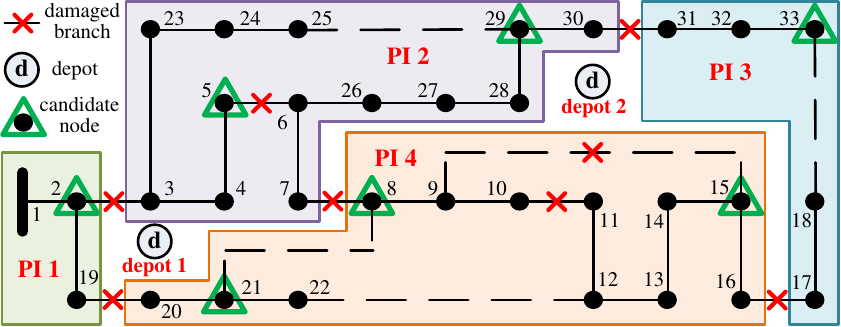}
  \vspace{-2.7mm}
  \caption{IEEE 33-node test system split into multiple PIs.}
  \vspace{-3.5pt}
  \label{33system}
\end{figure}

{\color{black}To coordinate with RC/MPS dispatch,
DS is dynamically reconfigured, e.g., to form microgrids.}
DS topology has to be radial in this process.
As researchers have extensively studied DS reconfigurution, the modeling of radiality
constraints is resolved to some extent. However, in this work we encounter a new
situation: The physical structure of the DS 
varies with the repair plan, and the source node distribution in the DS
varies with the allocation of MPSs. In general, we are reconfiguring a variable DS.
For example, a natural disaster split the DS in Fig.~\ref{33system} into
4 physical islands (PIs). Then, at some future time $t$ of the recovery process,
the number of PIs, the components of each PI, etc., are all variables that are not only {\color{black}dependent} on but
also to be co-optimized with RC and MPS dispatch decisions.
Existing methods to formulate radiality constraints in the literature are not applicable in this case.

\begin{figure}[t!]
  \centering
  \includegraphics[width=2.15in]{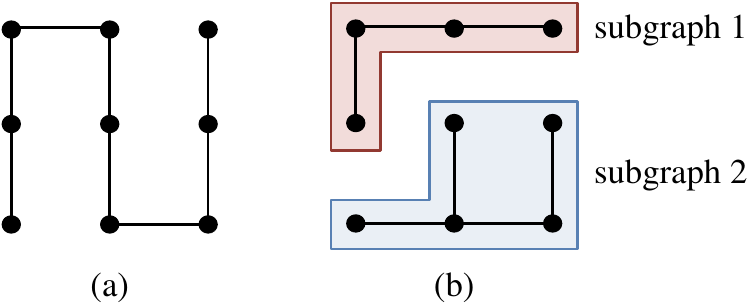}
  \vspace{-3.3mm}
  \caption{(a) A spanning tree; (b) A spanning forest.}
  \vspace{-8.5pt}
  \label{STandSF}
\end{figure}

Different from common DS reconfiguration problems that seek an optimal \emph{spanning tree}
(a radial topolgy connecting all nodes without loops),
in this work, with PIs, each feasible topology is a \emph{spanning foreast} (each subgraph being a spanning tree;
see Fig.~\ref{STandSF}).
Recent publications~\cite{Chen2016} and~\cite{Ding2017} on~microgrids formation are essentially constructing
a spanning foreast, too. However, their models are neither applicable here.
Also note that \cite{Ari2017} is essentially reconfiguring a variable DS. It uses the constraints from \cite{Jab2012}
to ensure radiality. However, a topology satisfying such constraints is not necessarily radial \cite{Ahm2015}.

Observing the difference and connection between spanning tree and spanning forest,
we propose the following method to express radiality constraints for a variable DS.
For the DS network, let $\lambda_{t}^{ij}$ be the connection status of branch $(i,j)$ at time $t$ (1 if closed, 0 if open);
let $u_{t}^{ij}$ be the operable status of branch $(i,j)$ at time $t$ (1 if can be used, 0 if damaged and unrepaired).
Introduce another fictitious network the same as the DS network but without damages.
Let $\boldsymbol{e_{t}}\triangleq\{e_{t}^{ij},\forall(i,j)\}$ be connection status of the fictitious network at time $t$.
Then, radiality constraints of our studied problem is expressed as:
\begin{equation}\label{radial_Constraints}
  \begin{aligned}
      (\ref{radial_Constraints}a):\ \boldsymbol{e_{t}}\in\boldsymbol{\Omega}, \forall t;
      \ \ \ \ (\ref{radial_Constraints}b):\ \lambda_{t}^{ij}\le e_{t}^{ij},\forall(i,j), \forall t;
      \\(\ref{radial_Constraints}c):\ \lambda_{t}^{ij}\le u_{t}^{ij},\forall(i,j), \forall t.
  \end{aligned}
\end{equation}
where $\boldsymbol{\Omega}$ is the set of spanning tree topologies of the fictitious network.
{\color{black} Note that (\ref{radial_Constraints}a) can be easily formulated by constraints based on
the single-commodity flow model \cite{Pop2002}, which
have been used extensively.
For space limit, we do not expand on the explicit formulation of (\ref{radial_Constraints}a),
and only write it symbolically here.}
In general, (\ref{radial_Constraints}a) requires $\boldsymbol{e_{t}}$ to form a fictitious spanning tree;
(\ref{radial_Constraints}b) restricts the DS to close only a subset of the branches in the fictitious spanning tree;
and (\ref{radial_Constraints}c) enforces inoperable branches to be open.
With \eqref{radial_Constraints}, $\boldsymbol{\lambda_{t}}\triangleq\{\lambda_{t}^{ij},\forall(i,j)\}$
forms a spanning forest in each time period.

{\color{black} The validity of~\eqref{radial_Constraints} is supported by the following property:
removing $l\ge0$ edges from a spanning tree leads to a spanning forest.
In other words, subgraphs of a spanning tree are also spanning trees, thus forming a spanning forest.}

Next, we model power dispatch of the DS.
Let $\delta_{t}^{i}$ be 1 if the load at node $i$ is restored at time $t$, 0 otherwise;
let $p_{t}^{i}$ and $q_{t}^{i}$ be the real and reactive power demand of node $i$ at time $t$, respectively;
let $P_{t}^{i}$ and $Q_{t}^{i}$ be the real and reactive power output of the MPS(s) at node $i$ at time $t$, respectively;
let $v_{t}^{i}$ be the squared voltage magnitude of node $i$ at time $t$;
let $\underline{v}^{i}$ and $\overline{v}^{i}$ be the minimum and maximum voltage value of node $i$, respectively;
let $pf_{t}^{ij}$ and $qf_{t}^{ij}$ be the real and reactive power flow on branch $(i,j)$ at time $t$, respectively;
let $r_{ij}$, $x_{ij}$ and $\overline{S}_{ij}$ be the resistance, reactance and apparent power capacity of branch $(i,j)$, respectively.
Then, we have:
\begin{equation}\label{DS0}
  \begin{aligned}
    \delta_{t}^{i} \le \delta_{t+1}^{i},
    \forall i, \forall t \le T-1.
  \end{aligned}
\end{equation}
\begin{equation}\label{DS1}
  \begin{aligned}
    P_{t}^{i} - \delta_{t}^{i} \cdot p_{t}^{i} + \sum_{(j,i) \in \boldsymbol{L}} pf_{t}^{ji} - \sum_{(i,j) \in \boldsymbol{L}} pf_{t}^{ij} = 0,
    \forall i, \forall t.
  \end{aligned}
\end{equation}
\begin{equation}\label{DS2}
  \begin{aligned}
    Q_{t}^{i} - \delta_{t}^{i} \cdot q_{t}^{i} + \sum_{(j,i) \in \boldsymbol{L}} qf_{t}^{ji} - \sum_{(i,j) \in \boldsymbol{L}} qf_{t}^{ij} = 0,
    \forall i, \forall t.
  \end{aligned}
\end{equation}
\begin{equation}\label{DS3}
  \begin{aligned}
    (\underline{v}^{i})^{2} \le v_{t}^{i} \le (\overline{v}^{i})^{2},
    \forall i, \forall t.
  \end{aligned}
\end{equation}
\begin{equation}\label{DS4}
  \begin{aligned}
    (pf_{t}^{ij})^{2} + (qf_{t}^{ij})^{2} \le \lambda_{t}^{ij} \cdot (\overline{S}_{ij})^{2},
    \forall (i,j), \forall t.
  \end{aligned}
\end{equation}
\begin{equation}\label{DS5}
  \begin{aligned}
    v_{t}^{i} - v_{t}^{j} \le (1 - \lambda_{t}^{ij}) \cdot K + 2 \cdot (r_{ij} \cdot pf_{t}^{ij} + x_{ij} \cdot qf_{t}^{ij}),
    \\\forall (i,j), \forall t.
  \end{aligned}
\end{equation}
\begin{equation}\label{DS6}
  \begin{aligned}
    v_{t}^{i} - v_{t}^{j} \ge (\lambda_{t}^{ij} - 1) \cdot K + 2 \cdot (r_{ij} \cdot pf_{t}^{ij} + x_{ij} \cdot qf_{t}^{ij}),
    \\\forall (i,j), \forall t.
  \end{aligned}
\end{equation}
where $\boldsymbol{L}$ is the set of DS branches; and $K$ is a large enough positive number.
As above, \eqref{DS0} prevent de-energizing loads that are already restored;
\eqref{DS1} and \eqref{DS2} require all nodes to satisfy real and reactive power balance conditions, respectively;
\eqref{DS3} specifies voltage magnitude limits;
\eqref{DS4} limit the apparent power on each branch by its capacity, and restrict both real and reactive power flow on a branch
to be zero if it is open;
\eqref{DS5} and \eqref{DS6} represent the power flow equation based on the DistFlow model \cite{Baran1989} \cite{Taylor2012}
for closed branches with $\lambda_{t}^{ij}=1$, and get relaxed for open branches with $\lambda_{t}^{ij}=0$.

Thus, with \eqref{radial_Constraints}-\eqref{DS6}, dynamic network reconfiguration
and power dispatch of the DS are coordinated to achieve service restoration
via strategies such as microgrid formation.

\section{The Co-Optimization Model}\label{CoOptModel}

{\color{black} DS service restoration is co-optimized with RC dispatch and MPS dispatch to enhance DS resilience.}
Let $\omega^{i}$ be the priority weight of the power demand at node $i$.
The objective function of the co-optimization is modeled as follows:
\begin{equation}\label{obj}
  \begin{aligned}
    \max \sum_{t} \left[ \sum_{i} \omega^{i} \cdot \delta_{t}^{i} \cdot q_{t}^{i}
    - \varepsilon \cdot \left( \sum_{k} \beta_{t}^{k} + \sum_{s} \beta_{t}^{s}
    \right) \right]
  \end{aligned}
\end{equation}
which maximizes the weighted sum of restored loads over~time,
and minimizes the total number of travels of RCs and MPSs.
Relative weights of the two objectives are adjusted by parameter $\varepsilon$,
which is set as a small value so that the first objective is still dominating.
We add the second objective for two main reasons:
1) To restrict the transportation of RCs and MPSs during time $t\sim T$ if all loads
are restored at some time $t<T$;
2) To select a dispatch strategy achieving the optimal recovery~effect~by a minimum number of travels of RCs and MPSs.

Next, we model the interdependence among RC dispatch, MPS dispatch and DS service restoration.
First, the relationship between RC dispatch and DS restoration is modeled.
Let branch $(i_{m},j_{m})$ denote the corresponding damaged component $m$;
let $\boldsymbol{L_{1}} \triangleq \{(i_{m},j_{m}), \forall m \in \boldsymbol{V_{1}}\}$
be the set of damaged branches;
let $\boldsymbol{L_{2}}$ be the set of branches with switches. Then, we have:
\begin{equation}\label{dependence1}
  \begin{aligned}
    u_{t}^{ij} = 1,
    \forall (i,j) \in \boldsymbol{L} \setminus \boldsymbol{L_{1}},
    \forall t.
  \end{aligned}
\end{equation}
\begin{equation}\label{dependence2}
  \begin{aligned}
    u_{t+1}^{i_{m}j_{m}} \le \sum_{k}z_{m,t}^{k},
    \forall m \in \boldsymbol{V_{1}},
    \forall t \le T-1.
  \end{aligned}
\end{equation}
\begin{equation}\label{dependence3}
  \begin{aligned}
    \lambda_{t}^{ij} = 1,
    \forall (i,j) \in \boldsymbol{L} \setminus \{\boldsymbol{L_{1}} \cup \boldsymbol{L_{2}}\},
    \forall t.
  \end{aligned}
\end{equation}
As above, \eqref{dependence1} states that an intact branch is always operable;
\eqref{dependence2} indicates that a damaged branch is operable only if it is repaired
by one of the RCs in the previous time period;
\eqref{dependence3} enforces undamaged branches without switches to be closed.

Second, the relationship between MPS dispatch and DS~service restoration is
incorporated by the following equations:
\begin{equation}\label{dependence4}
  \begin{aligned}
    P_{t}^{i} = \sum_{s \in \boldsymbol{M_{1}}} \alpha_{i,t}^{s} \cdot gp_{t}^{s}
    + \sum_{s \in \boldsymbol{M_{2}}} \alpha_{i,t}^{s} \cdot (dp_{t}^{s} - cp_{t}^{s}),
    \\\forall i \in \boldsymbol{N'}, \forall t.
  \end{aligned}
\end{equation}
\begin{equation}\label{dependence5}
  \begin{aligned}
    Q_{t}^{i} = \sum_{s \in \boldsymbol{M}} \alpha_{i,t}^{s} \cdot gq_{t}^{s},
    \forall i \in \boldsymbol{N'}, \forall t.
  \end{aligned}
\end{equation}
\begin{equation}\label{dependence6}
  \begin{aligned}
    P_{t}^{i} = Q_{t}^{i} = 0,
    \forall i \in \boldsymbol{N} \setminus \boldsymbol{N'}, \forall t.
  \end{aligned}
\end{equation}
With \eqref{dependence4} and \eqref{dependence5}, $P_{t}^{i}$ and $Q_{t}^{i}$
are derived by summing the real and reactive power outputs of MPSs connected to node $i$ at time $t$, respectively;
\eqref{dependence6} enforces $P_{t}^{i}$ and $Q_{t}^{i}$ to be zero for nodes that are not for MPS connection.

Note that although \eqref{dependence1}-\eqref{dependence6} do not show direct
relationship between RC dispatch and MPS dispatch, they are intrinsically interrelated
in the studied problem.
They need to be coordinated, i.e., co-optimized,
to attain better service restoration of the DS.

Thus, we arrive at the co-optimization model as follows:
\begin{equation*}
  \begin{aligned}
    Objective:\ \ & \eqref{obj}; \\
    Constraints:\ \ & \eqref{AtMostOneRepair},\eqref{AuxRC},\eqref{TravelTimeTrans},
    \eqref{Repair1}{-}\eqref{DS6},\eqref{dependence1}{-}\eqref{dependence6}; \\
    Variables:\ \ & a_{m,t}^{k}, \beta_{t}^{k}, z_{m,t}^{k};\\
    & \alpha_{i,t}^{s}, \beta_{t}^{s}, gp_{t}^{s}, gq_{t}^{s}, c_{t}^{s}, d_{t}^{s}, cp_{t}^{s}, dp_{t}^{s}, soc_{t}^{s};\\
    & \delta_{t}^{i}, P_{t}^{i}, Q_{t}^{i}, v_{t}^{i}, \lambda_{t}^{ij}, u_{t}^{ij}, e_{t}^{ij}, pf_{t}^{ij}, qf_{t}^{ij}.
  \end{aligned}
\end{equation*}
Above, the list of variables are also provided for clarity.

\section{Solution Method}\label{Algorithms}

\subsection{Linearization Techniques}

The proposed co-optimization model is a non-convex MINLP, as \eqref{dependence4} and \eqref{dependence5} have non-linear and non-convex terms
such~as $\alpha_{i,t}^{s} \cdot gp_{t}^{s}$.
They are linearized by the McCormick envelopes \cite{McCormick1976}.
As $\alpha_{i,t}^{s}$ is binary, and the involved continuous variables have explicit lower and upper bounds,
the linearization is equivalent. By doing this, the co-optimization model
is also convexified into a MISOCP.
We further linearize \eqref{DS4} using the technique~in \cite{ChenXin2016}.
The co-optimization model thus becomes a MILP, which can be solved
by off-the-shelf solvers such as Gurobi.
As the involved linearizations are straightforward,
we do not elaborate on the detailed reformulations for space limit.

\vspace{-2pt}
\subsection{Pre-Assigning a Minimal Set of Repair Tasks to Depots}

To reduce the MILP co-optimization model's computational complexity,
we follow \cite{Ari2017} to cluster and pre-assign repair tasks to depots.
Let $\psi_{n}^{m}$ be 1 if damaged component $m$ is assigned to depot $n$, 0 otherwise;
let $l_{mn}$ be the distance between damaged component $m$ and depot $n$;
let $\boldsymbol{\Phi_{n}}$ be the set of RCs in depot $n$.
Then, the pre-assignment of repair tasks to depots is determined by the small MILP model as follows:
\begin{equation}\label{PA1}
  \begin{aligned}
    \min_{\psi_{n}^{m},u^{ij},\lambda^{ij},e^{ij},\atop pf^{ij},qf^{ij},v^{i}.} \sum_{m \in \boldsymbol{V_{1}}} \sum_{n \in \boldsymbol{V_{2}}}
    \psi_{n}^{m} \cdot l_{mn}
  \end{aligned}
\end{equation}
\begin{equation}\label{PA2}
  \begin{aligned}
    s.t. \ \ \ \ \sum_{m \in \boldsymbol{V_{1}}} \psi_{n}^{m} \cdot rs_{m} \le \sum_{k \in \boldsymbol{\Phi_{n}} } RS^{k},
    \forall n \in \boldsymbol{V_{2}}. \ \ \
  \end{aligned}
\end{equation}
\begin{equation}\label{PA3}
  \begin{aligned}
    \sum_{n \in \boldsymbol{V_{2}}} \psi_{n}^{m} \le 1, \forall m \in \boldsymbol{V_{1}}.
  \end{aligned}
\end{equation}
\begin{equation}\label{PA4}
  \begin{aligned}
    u^{i_{m}j_{m}} \le \sum_{n \in \boldsymbol{V_{2}}} \psi_{n}^{m}, \forall m \in \boldsymbol{V_{1}}.
  \end{aligned}
\end{equation}
\begin{equation}\label{PA5}
  \begin{aligned}
    \sum_{(j,i) \in \boldsymbol{L}} pf^{ji} - \sum_{(i,j) \in \boldsymbol{L}} pf^{ij} - p^{i} = 0,
    \forall i.
  \end{aligned}
\end{equation}
\begin{equation}\label{PA6}
  \begin{aligned}
    \sum_{(j,i) \in \boldsymbol{L}} qf^{ji} - \sum_{(i,j) \in \boldsymbol{L}} qf^{ij} - q^{i} = 0,
    \forall i.
  \end{aligned}
\end{equation}
\begin{equation}\label{PA7}
    \eqref{radial_Constraints},\eqref{DS3}{-}\eqref{DS6},\eqref{dependence1},\eqref{dependence3} \
    (\mathrm{without \ time \ subscript}\ t).
\end{equation}
Specifically, \eqref{PA1} minimizes the sum of distances between~damaged components
and their assigned depots, so as to reduce the transportation time of RCs and thus enhance DS restoraiton;
\eqref{PA2} requires each depot to have adequate resources
to complete its assigned repair tasks;
and \eqref{PA3} states that each damaged component is assigned to at
most one depot.
Note that, different from \cite{Ari2017}, here we do not demand all tasks to be assigned.
The above model just seeks and assigns \emph{a minimal set of repair tasks} that can fully restore all DS loads
without MPSs. Thus, \eqref{PA4}-\eqref{PA7} are added:
\eqref{PA4} states that a damaged branch will be operable only if it is assigned to be repaired;
\eqref{PA5} and \eqref{PA6} represent the power balance conditions to fully supply all loads without MPSs;
\eqref{PA7} are the other DS operational constraints.
Parameters $p^{i}/q^{i}$ can be set as the maximum power demands of DS nodes. And note that power outputs of the substation
are implicitly modeled as $pf^{i_{g}j}/qf^{i_{g}j}$, where $i_{g}$ denotes the substation node.
Though without the time subscript $t$,
other involved variables and constraints are introduced hereinbefore.
By selecting and assigning a minimal set of repair tasks rather than all tasks to depots, computational complexity of
the co-optimization model
is further reduced.
Another advantage is that,
after all loads are restored
by a minimal set of repair~tasks, RCs can be better scheduled to help recover other DSs.

\vspace{-2pt}
\subsection{Selecting Candidate Nodes for MPS Connection}

The co-optimization model's computation time also partially depends on the
number of candidate nodes for MPS connection in the DS.
Here, one candidate node is selected from each PI.

First, we need to find the set of candidate nodes in~each~PI.
Let $\boldsymbol{I}$ be the set of all candidate nodes,
i.e., nodes meeting~the requirements of MPS connection, e.g., facility requirements~\cite{Lei2016};
let $\boldsymbol{I_{\sigma}}$ be the set of candidate nodes in PI $\sigma$;
let $v2s$ be an~operation transforming a vector into a set, e.g., $[0,1,2] \rightarrow \{0,1,2\}$;
let $\odot$ be the symbol for Hadamard product, i.e., element-wise multiplication;
let $\boldsymbol{\Theta}$ be the adjacency matrix
($\boldsymbol{\Theta_{ij}}=1$ if there is an intact branch $(i,j)$, whether
closed or open; $\boldsymbol{\Theta_{ij}}=0$ otherwise);
let $\boldsymbol{\Theta_{i}}$ be the $i$th row of $\boldsymbol{\Theta}$.
Algorithm~\ref{alg1} as follows is proposed to find $\boldsymbol{I_{\sigma}}$ for each PI $\sigma$:
\begin{algorithm}
\caption{Find the set of candidate nodes in each PI}\label{alg1}
Input: $\boldsymbol{N}$, $\boldsymbol{\Theta}$, $\boldsymbol{I}$; \
Output: $\boldsymbol{I_{\sigma}},\forall \sigma.$
\begin{algorithmic}[1]
\State $\sigma \gets 1$, $\boldsymbol{\mathcal{N}} \gets \emptyset$;
\While {$\boldsymbol{N} \setminus \boldsymbol{\mathcal{N}} \neq \emptyset$}
  \State randomly select  $i \in \boldsymbol{N} \setminus \boldsymbol{\mathcal{N}}$,
          $\boldsymbol{I_{\sigma}} \gets \{i\}$,
          $\boldsymbol{\mathcal{I}_{\sigma}} \gets \{i\}$;
  \State $\boldsymbol{C_{i}} \gets v2s(\boldsymbol{\Theta_{i}} \odot [1,...,N])$,
          $\boldsymbol{I_{\sigma}} \gets \{\boldsymbol{I_{\sigma}} \cup \boldsymbol{C_{i}}\} \setminus \{0\}$;
  \If {$\boldsymbol{I_{\sigma}} \setminus \boldsymbol{\mathcal{I}_{\sigma}} \neq \emptyset$}
    \State randomly select $i \in \boldsymbol{I_{\sigma}} \setminus \boldsymbol{\mathcal{I}_{\sigma}}$,
            $\boldsymbol{\mathcal{I}_{\sigma}} \gets \boldsymbol{\mathcal{I}_{\sigma}} \cup \{i\}$;
    \State  \textbf{go to} \emph{step 4};
  \EndIf
  \State $\boldsymbol{\mathcal{N}} \gets \boldsymbol{\mathcal{N}} \cup \boldsymbol{I_{\sigma}}$,
          $\boldsymbol{I_{\sigma}} \gets \boldsymbol{I_{\sigma}} \cap \boldsymbol{I}$,
          $\sigma \gets \sigma + 1$;
   \State \textbf{go to} \emph{step 3};
\EndWhile
\end{algorithmic}
\end{algorithm}

Then, the final selection of candidate nodes is determined by the
small integer programming (IP) model as follows:
\begin{equation}\label{CN1}
  \begin{aligned}
    \min_{\phi_{i}} \frac{1}{2}
    \sum_{\sigma}
    \sum_{i \in \boldsymbol{I_{\sigma}}}
    \sum_{j \in \boldsymbol{I} \setminus \boldsymbol{I_{\sigma}}}
    \phi_{i} \cdot \phi_{j} \cdot l_{ij}
  \end{aligned}
\end{equation}
\begin{equation}\label{CN2}
  \begin{aligned}
    s.t. \ \ \
    \sum_{i \in \boldsymbol{I_{\sigma}}}
    \phi_{i} = 1,
    \forall \sigma \ with \ I_{\sigma} \neq 0.
    \ \ \ \ \ \ \
  \end{aligned}
\end{equation}
where $\phi_{i}=1$ if candidate node $i$ is selected, $\phi_{i}=0$ otherwise;
$l_{ij}$ is the distance between nodes $i$ and $j$;
$I_{\sigma}$ is the cardinality of $\boldsymbol{I_{\sigma}}$;
\eqref{CN1} minimizes the sum of distances between selected candidate
nodes, in order to reduce the transportation time of MPSs and thus enhance DS restoration;
\eqref{CN2} enforces one node to be selected from each PI that has at least one candidate~node.
Again, by the McCormick envelopes \cite{McCormick1976}, terms $\phi_{i} \cdot \phi_{j}$ can~be equivalently linearized.

\section{Case Studies}\label{cases}

In this section, the proposed co-optimization method for~disaster recovery
logistics is demonstrated on two systems. We~use a computer with an Intel
i5-4278U processor and 8GB memory. Involved MILP and IP problems are solved
by Gurobi 7.5.2.

\subsection{Case I: IEEE 33-Node Test System}

For this DS, we consider a scenario with 8 branches~damaged by the
natural disaster (see Fig.~\ref{33system}). The system has 2 depots~and 2 RCs, i.e., RC 1 in depot 1
and RC 2 in depot 2. Resource capacity of RCs is set as 8.
Both the number of resources and time periods required to repair
different~damaged~components~vary from 1 to 4.
The DS also has 2~MPSs,~i.e.,~500~kW/400 kVar MEG 1, and 300 kW/300 kWh MESS 1.
Priority weights of loads and the travel time/distance data are randomly generated.
Branches 9-10, 9-15, 12-22, 14-15, 18-33 25-29, 28-29 and 30-31 are equipped with remote-contorlled switches \cite{Lei2017}.
Some other data can be found in \cite{Baran1989}.

First, a minimal set of repair tasks is sought and pre-assigned to depots.
The pre-assignment model \eqref{PA1}-\eqref{PA7} is solved within 0.13 s.
Table~\ref{PreA} lists the results.
In virtue of the back-up branches, i.e., the normally open branches,
repairing 4 damaged branches is sufficient to fully restore all loads.
The repair of the other 4 un-assigned damaged branches can be arranged later for
an objective other than maximizing the sum of restored loads.

Second, candidate nodes for MPS connection are selected.
Let $\boldsymbol{N'}=\{2,5,8,15,21,29,33\}$ be the set of all candidate nodes.
Algorithm~\ref{alg1}, which takes 0.14 s, accurately detects~that the system is
split into 4 PIs containing candidate nodes $\{2\}$, $\{5,29\}$,
$\{8,15,21\}$ and $\{33\}$, respectively.
Note that although these results are obvious in this case,
Algorithm~\ref{alg1} is necessary to automatically generate such data for both small and large systems.
Then, the candidate node selection problem \eqref{CN1}-\eqref{CN2} is solved in
0.04 s. Nodes $\{2,8,29,33\}$ are selected.

\begin{figure}[t!]
    \centering
    \includegraphics[width=2.4in]{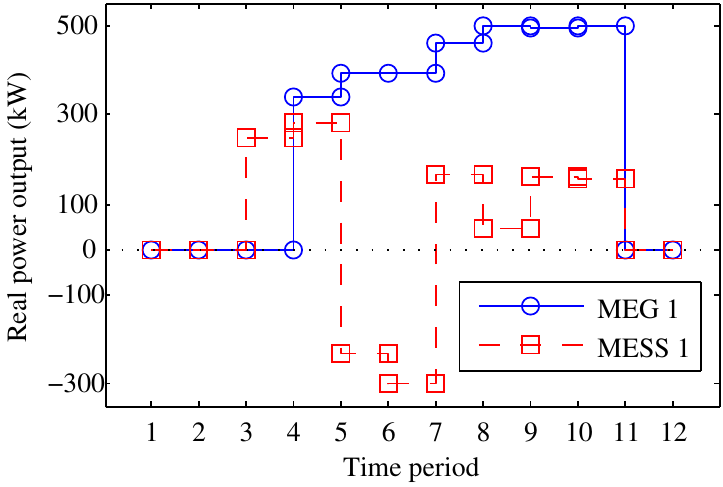}
    \vspace{-2.75mm}
    \caption{Real power outputs of MPSs in each time period  (Case I).}
    \label{MPSpower}
\end{figure}

Third, after pre-assigning repair tasks and selecting candidate nodes, we solve the co-optimization
problem. Let $T=12$~and $\Delta t = 0.5\ \mathrm{hr}$; let $vol_{i}=2$ for all candidate nodes;
let $\alpha^{s}_{2,1}=1$ for both MPSs. The co-optimization model is~solved~in~2.26~s.
Table~\ref{RCRS}~and~\ref{RSMPS} list the dispatch of RCs and MPSs, respectively.
Symbols ``$\rightarrow$'' or ``\text{\sffamily x}'' mean that this RC/MPS
is being~transported or has stopped working, respectively.
MPSs change their locations only once for this small system.
We will see more dynamic dispatch of MPSs in Case II.
Fig.~\ref{MPSpower} depicts MPSs' real power outputs.
Table~\ref{RCSactions} lists switch actions of the DS.

\begin{figure}[t!]
    \centering
    \includegraphics[width=3.4in]{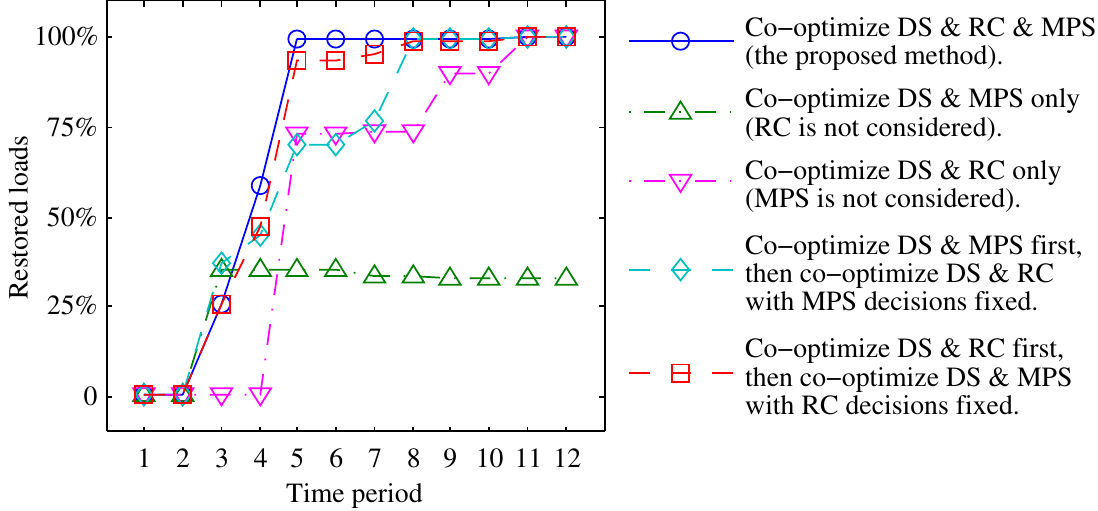}
    \vspace{-2.75mm}
    \caption{Restored loads over time for different logistics strategies (Case I).}
    \label{loadcurves}
\end{figure}

The proposed co-optimization method for disaster recovery logistics is also
compared to other logistics strategies.
As~shown in Fig~\ref{loadcurves}, the proposed method co-optimizing both RC~dispatch
and MPS dispatch with DS restoration has the best performance.

\begin{table}[t!]
  \centering
  \caption{Pre-Assignment of A Minimal Set of Repair Tasks (Case I)}
  \vspace{-2.5mm}
    \begin{tabular}{c|c|c}
    \hline 
        ------  & Depot 1 & Depot 2 \\
    \hline 
    Pre-assigned repair tasks & branches 2-3, 5-6 & branches 7-8, 30-31 \\
    \hline 
    \hline 
    Un-assigned repair tasks & \multicolumn{2}{|c}{banches 10-11, 16-17, 19-20, 9-15} \\
    \hline 
    \end{tabular}%
  \label{PreA}%
\end{table}%

\begin{table}[t!]
  \centering
  \caption{Routing and Scheduling of RCs (Case I)}
  \vspace{-2.5mm}
    \setlength{\tabcolsep}{2.9pt}
    \begin{tabular}{c|c|c|c|c|c|c|c}
    \hline
    \multirow{2}{*}{RC 1} & Time period & 0     & 1     & 2$\sim$4   & 5$\sim$7   & 8$\sim$10  & 11$\sim$12 \\
    \cline{2-8}
          & Dispatch & depot 1 & $\rightarrow$       & branch 2-3 &   $\rightarrow$     & branch 5-6 &  \text{\sffamily x} \\
    \hline
    \hline
    \multirow{2}{*}{RC 2} & Time period & 0     & 1     & 2$\sim$3   & 4$\sim$5   & 6$\sim$7   & 8$\sim$12 \\
    \cline{2-8}
          & Dispatch & depot 2 &   $\rightarrow$     & branch 7-8 &    $\rightarrow$    & branch 30-31 &  \text{\sffamily x} \\
    \hline
    \end{tabular}%
  \label{RCRS}%
\end{table}%

\begin{table}[t!]
  \centering
  \caption{Routing and Scheduling of MPSs (Case I)}
  \vspace{-2.5mm}
    \begin{tabular}{c|c|c|c|c|c}
    \hline
    \multirow{2}[0]{*}{MEG 1} & Time period & 0     & 1$\sim$3   & 4$\sim$10 & 11$\sim$12 \\ \cline{2-6}
              & Dispatch & node 2 &   $\rightarrow$    & node 33 & \text{\sffamily x} \\
    \hline
    \hline
    \multirow{2}[0]{*}{MESS 1} & Time period & 0     & 1$\sim$2   & 3$\sim$10 & 11$\sim$12 \\ \cline{2-6}
          & Dispatch & node 2 &   $\rightarrow$    & node 29 & \text{\sffamily x} \\
    \hline
    \end{tabular}%
  \label{RSMPS}%
\end{table}%

\begin{table}[t!]
  \centering
  \caption{Dynamic Network Reconfiguration of the DS (Case I)}
  \vspace{-2.5mm}
    \begin{tabular}{c|c|c|c}
    \hline
    Time period & 3     & 4     & 11 \\
    \hline
    Switch actions & close 25-29 & close 8-21, 12-22, 18-33 & open 28-29 \\
    \hline
    \end{tabular}%
  \label{RCSactions}%
\end{table}%

\begin{figure*}[!t]
    \centering
    \includegraphics[width=7.146in]{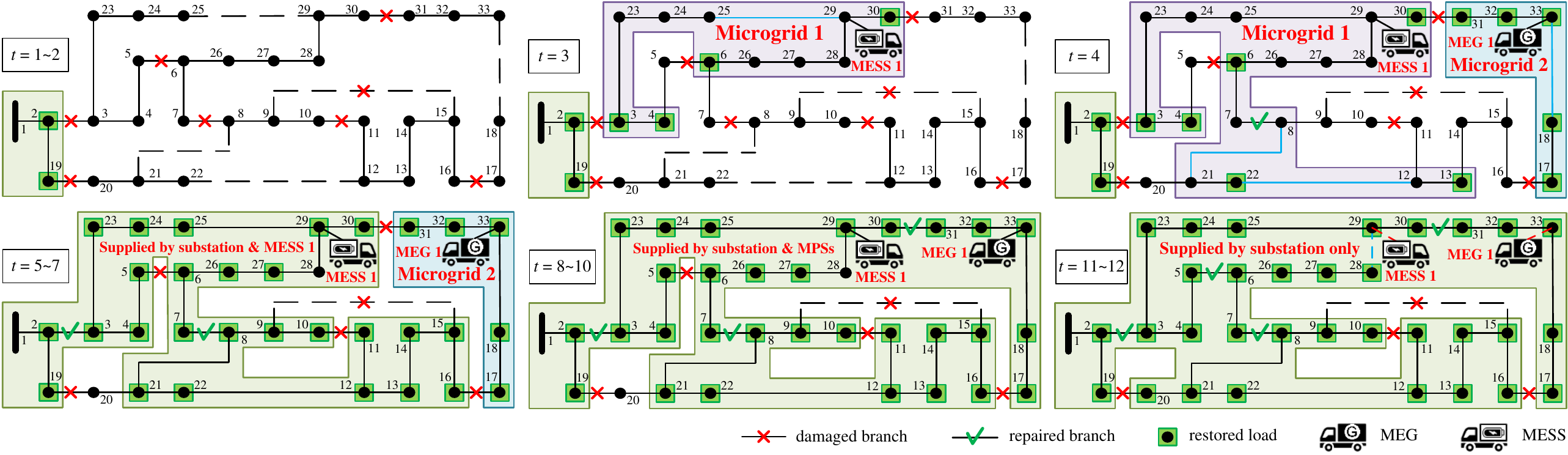} 
    \vspace*{-7.5mm}
    \caption{DS service restoration process co-optimized with both RC dispatch and MPS dispatch (Case I).}
    \label{restorationProcess}
\end{figure*}

\begin{table*}[htbp]
  \centering
  \caption{Routing and Scheduling of RCs (Case II)}
  \vspace{-2.5mm}
    \begin{tabular}{c|c|c|c|c|c|c|c|c|c|c|c}
    \hline
    \multirow{2}[0]{*}{RC 1} & Time period & 0     & 1     & 2$\sim$4 & 5$\sim$6 & 7$\sim$9 & 10$\sim$11 & 12$\sim$13 & 14    & 15    & 16 \\
    \cline{2-12}
              & Dispatch & depot 1 & $\rightarrow$ & branch 9-15 & $\rightarrow$ & branch 72-73 & $\rightarrow$ & branch 88-89 & $\rightarrow$ & branch 106-107 & \text{\sffamily x} \\
    \hline
    \hline
    \multirow{2}[0]{*}{RC 2} & Time period & 0     & 1     & 2$\sim$3 & 4     & 5     & 6     & 7     & 8$\sim$9 & 10$\sim$11 & 12$\sim$16 \\
    \cline{2-12}
          & Dispatch & depot 1 & $\rightarrow$ & branch 103-108 & $\rightarrow$ & branch 39-41 & $\rightarrow$ & branch 15-56 & $\rightarrow$ & branch 65-72 & \text{\sffamily x} \\
    \hline
    \hline
    \multirow{2}[0]{*}{RC 3} & Time period & 0     & 1     & 2     & 3     & 4$\sim$5 & 6     & 7$\sim$8 & 9     & 10$\sim$11 & 12$\sim$16 \\
    \cline{2-12}
          & Dispatch & depot 2 & $\rightarrow$ & branch 27-33 & $\rightarrow$ & branch 92-93 & $\rightarrow$ & branch 99-101 & $\rightarrow$ & branch 17-19 & \text{\sffamily x} \\
    \hline
    \hline
    \multirow{2}[0]{*}{RC 4} & Time period & 0     & 1     & 2$\sim$3 & 4     & 5     & 6     & 7$\sim$8 & 9$\sim$10 & 11$\sim$12 & 13$\sim$16 \\
    \cline{2-12}
          & Dispatch & depot 2 & $\rightarrow$ & branch 43-45 & $\rightarrow$ & branch 74-75 & $\rightarrow$ & branch 25-26 & $\rightarrow$ & branch 113-116 & \text{\sffamily x} \\
    \hline
    \end{tabular}%
  \label{RC123}%
\end{table*}%

The recovery process is decribed in details by Fig.~\ref{restorationProcess}.
At $t=3$, MESS 1 arrives at and is connected to node 29. Closing branch 25-29,
microgrid 1 is formed to restore nodes 3, 4, 6 and 30.
At $t=4$, with branch 7-8 repaired, microgrid~1 can further restore nodes 13 and 22
by closing branches 8-21 and 12-22. At $t=4$, MEG 1 also arrives at and is connected to node~33.
Closing branch 18-33, microgrid 2 is formed to restore all loads in it.
At $t=5$, with branch 2-3 repaired, microgrid~1 gets connected to the substation.
Thus, MESS 1 can~get~charged~at $t=5\sim6$ to supply future peak loads.
At $t=8$, with branch 30-31 repaired, microgrid 2 also gets connected to the main grid.
However, due to operational constraints (low voltage of node 16 and large power flow on branch 24-25),
nodes 20 and 28 cannot be restored untile branch 5-6 is repaired at $t=11$,
when branch 28-29 is opened to avoid loop,
and both MPSs are disconnected from the DS as they are no longer necessary.

\begin{figure}[t!]
  \centering
  \includegraphics[width=3.3in]{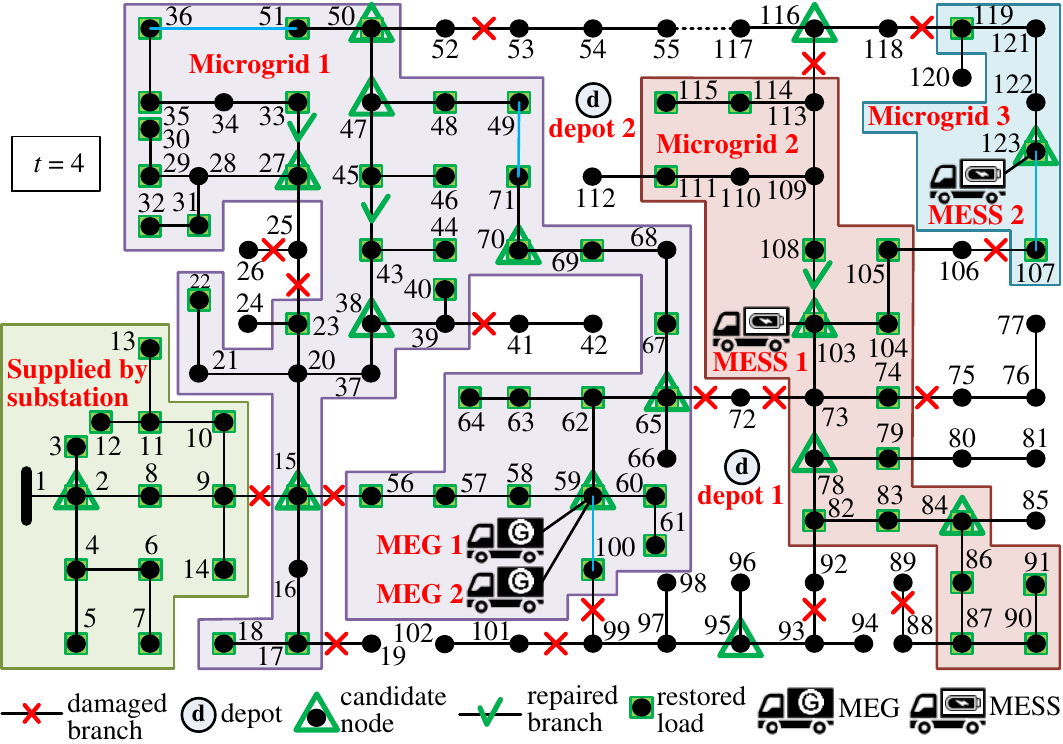}
  \vspace{-2.5mm}
  \caption{IEEE 123-node test system and its restoration at $t=4$.}
  \label{123system}
\end{figure}

\begin{table}[t!]
  \centering
  \caption{Routing and Scheduling of MPSs (Case II)}
  \vspace{-2.5mm}
  \setlength{\tabcolsep}{1.625pt}
    \begin{tabular}{c|c|c|c|c|c|c|c}
    \hline
    \multirow{2}[0]{*}{MEG 1} & Time period & 0     & 1$\sim$2   & 3$\sim$4   & 5     & 6$\sim$11  & 12$\sim$16 \\
    \cline{2-8}
              & Dispatch & node 2 & $\rightarrow$ & node 59 & $\rightarrow$ & node 103 & \text{\sffamily x} \\
    \hline
    \hline
    \multirow{4}[0]{*}{MEG 2} & Time period & 0     & 1$\sim$2   & 3$\sim$4   & 5   & 6$\sim$11  & 12$\sim$13    \\
    \cline{2-8}
          & Dispatch & node 2 & $\rightarrow$ & node 59 & $\rightarrow$ &   node 103 & $\rightarrow$  \\
          \cline{2-8}
          & Time period  & 14$\sim$15 & 16    &     \multicolumn{4}{c}{\multirow{2}{*}{---------}}   \\
          \cline{2-4}
          & Dispatch &  node 123 & \text{\sffamily x}     &      \multicolumn{4}{c}{}  \\
          \hline
          \hline
    \multirow{4}[0]{*}{MESS 1} & Time period & 0$\sim$1     & 2$\sim$3   & 4$\sim$5   & 6$\sim$7   & 8     & 9\\
    \cline{2-8}
          & Dispatch & node 123 & $\rightarrow$ & node 103 & $\rightarrow$ & node 38 & $\rightarrow$ \\
          \cline{2-8}
          & Time period    & 10$\sim$12 & 13$\sim$16 &       \multicolumn{4}{c}{\multirow{2}{*}{---------}}    \\
          \cline{2-4}
          & Dispatch & node 116 & \text{\sffamily x}   &       \multicolumn{4}{c}{}   \\
          \hline
          \hline
    \multirow{2}[0]{*}{MESS 2} & Time period & 0$\sim$14  & 15$\sim$16 &      \multicolumn{4}{c}{\multirow{2}{*}{---------}}  \\
    \cline{2-4}
          & Dispatch & node 123 & \text{\sffamily x}   & \multicolumn{4}{c}{}    \\
    \hline
    \end{tabular}%
  \label{MPS123}%
\end{table}%

\begin{figure}[t!]
    \centering
    \includegraphics[width=2.4in]{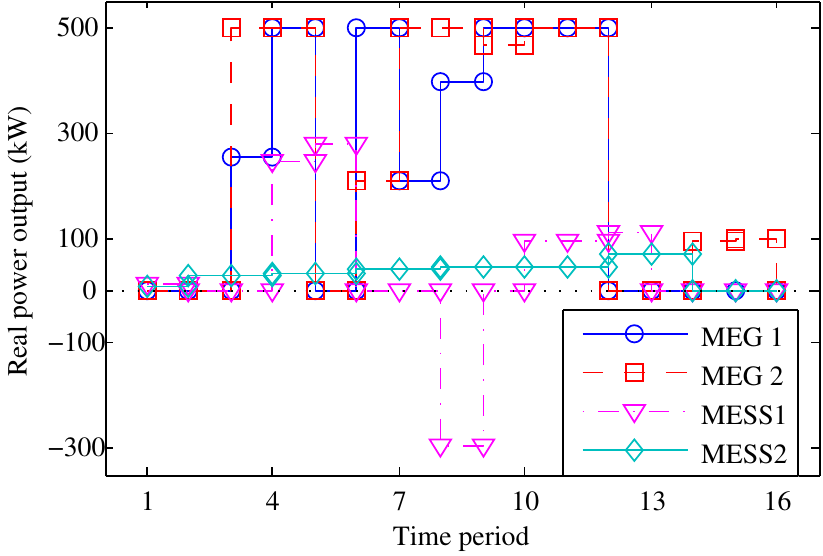}
    \vspace{-2.75mm}
    \caption{Real power outputs of MPSs in each time period  (Case II).}
    \label{MPSpower123}
\end{figure}

\subsection{Case II: IEEE 123-Node Test System}

A scenario with 20 branches damaged by the natural disaster
is considered for the second test system in Fig.~\ref{123system}.
It has 2~depots (each with 2 RCs), 4 MPSs (2 MEGs and 2 MESSs),
and 12 branches equipped with remote-controlled switches \cite{Lei2017}.
Let $T=16$.
Some other data can be found in Case I or in \cite{123system}.

\begin{figure}[t!]
    \centering
    \includegraphics[width=3.4in]{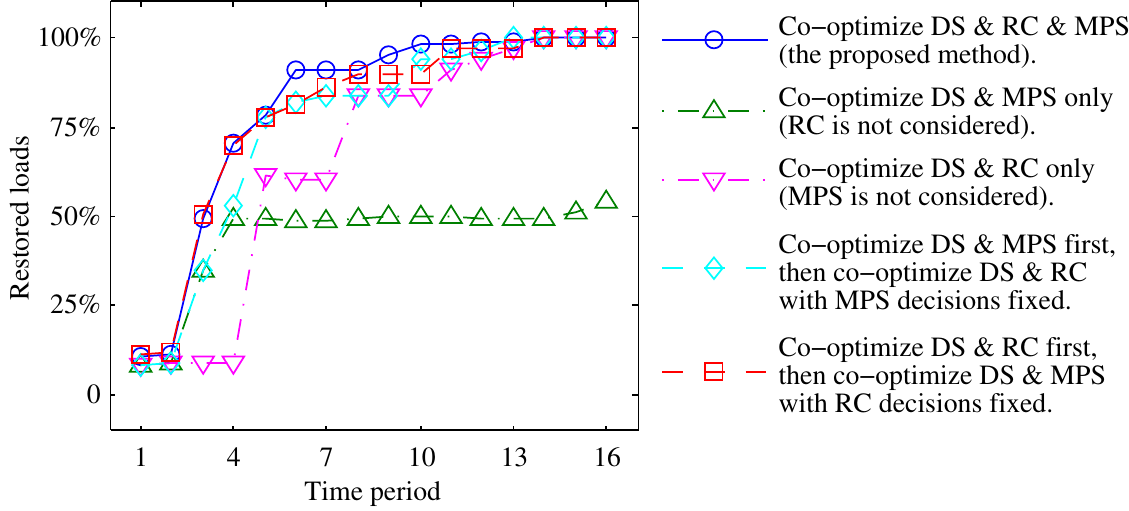}
    \vspace{-2.75mm}
    \caption{Restored loads over time for different logistics strategies (Case II).}
    \label{loadcurves123}
\end{figure}

First, the pre-assignment model \eqref{PA1}-\eqref{PA7} is solved in 0.43~s.
Next, Algorithm 1 is run for 0.45 s,
and the candidate node selection problem \eqref{CN1}-\eqref{CN2} is solved in 0.28 s.
It selects nodes $\{2,27,38,59,95,103,116,123\}$.
At last, the co-optimization problem is solved in 1183 s.
For space limit, we do not elaborate on the recovery process.
Only the DS's restoration progress~at $t=4$ is depicted in Fig~\ref{123system} for illustration.
Table~\ref{RC123} and \ref{MPS123} list the dispatch of RCs and MPSs, respectively.
Table~\ref{RC123} implicitly includes the results for the pre-assignment of a minimal set of repair tasks, too.
Fig.~\ref{MPSpower123} depicts real power outputs of MPSs.
Table~\ref{RCSactions123} lists the DS's switch actions in its dynamic network reconfiguration.
Fig.~\ref{loadcurves123} again demonstrates the effectiveness and superiority of
the proposed co-optimization method for disaster recovery logistics.
The tables and figures indicate that, the proposed method is especially effective in coordinating RC dispatch and MPS dispatch
to restore loads by dynamically forming microgrids in the DS.
The microgrids are powered by MPSs, and reconfigured and extended by switch actions of the DS
and repair actions of RCs.

\begin{table}[t!]
\centering
  \caption{Dynamic Network Reconfiguration of the DS (Case II)}
  \vspace{-2.5mm}
  \label{RCSactions123}%
      \setlength{\tabcolsep}{4pt}
\begin{tabular}{c|c|c|c|c}
  \hline
      Time period              & 1 &      3             &      8 & 10             \\
 \hline
 Switch & \multirow{2}{*}{close 107-123} & close 36-51, & \multirow{2}{*}{open 49-71} & \multirow{2}{*}{close 55-117} \\
 actions &                   & 49-71, 59-100  &                   &             \\
 \hline
\end{tabular}
\end{table}

\subsection{Computational Efficiency}
Table~\ref{computationTime1} reports the computation time of different solution
methods for the co-optimization model.
For Case I, all methods attain the same solution.
That is, preprocessing of the proposed method does not result in sub-optimality.
For Case II, only the proposed method solves the co-optimization problem in 2 hr.
Specifically, it is indicated that, compared to pre-assigning all repair tasks,
pre-assigning a minimal set of repair tasks
can better improve the computational efficiency.

\begin{table}[t!]
  \centering
  \caption{Computation Time of the Co-Optimization Problem (Case I \& II)}
  \vspace{-2.5mm}
        \setlength{\tabcolsep}{4.5pt}
    \begin{tabular}{c|c|c}
    \hline
    Solution methods & Case I & Case II \\
    \hline
    Proposed method & 2.26 s & 1183 s \\
    \hline
    Without any preprocessing & 66.20 s & gap = 97.20\% at 2 hr \\
    \hline
    Only with $\boldsymbol{N'}$ reduced & 37.59 s &  gap = 59.92\% at 2 hr \\
    \hline
    Only with all repairs pre-assigned & 17.69 s & gap = 32.60\% at 2 hr \\
    \hline
    Only with minimal repairs pre-assigned & 6.67 s &  gap = 5.34\% at 2 hr \\
    \hline
    \end{tabular}%
  \label{computationTime1}%
\end{table}%

\section{Conclusion}\label{conclusion}

To enhance DS resilience,
this paper builds a co-optimization method for~disaster recovery logistics.
RC dispatch, MPS~dispatch and DS operation are jointly cooperated for electric~service restoration.
A MINLP model is formulated to~attain~resilient strategies that involve the
routing and scheduling of~ RCs and MPSs, dynamic network reconfiguration of the DS, etc.~The model is transformed into a MILP,
and preprocessed~to~reduce computational complexity.
Case studies demonstrate that the co-optimization method for~disaster recovery logistics efficiently
improves DS service restoration,
{\color{black} especially by dynamic formation of microgrids that
are powered by MPSs and topologized by repair actions of RCs
and reconfiguration of the DS.}


%

%



\ifCLASSOPTIONcaptionsoff
  \newpage
\fi



%
\bibliographystyle{IEEEtran}
\bibliography{IEEEabrv,references}

%








\end{document}